\documentclass[12pt]{amsproc}%
\usepackage{amsfonts}
\usepackage{amsmath}
\usepackage{amssymb}
\usepackage{graphicx}%
\renewcommand {\theequation} {\@arabic\c@equation}
 \theoremstyle{plain}

\newtheorem{corollary}{Corollary}

\newtheorem{definition}{Definition}

\newtheorem{lemma}{Lemma}

\newtheorem{proposition}{Proposition}
\newtheorem{remark}{Remark}

\newtheorem{theorem}{Theorem}
\numberwithin{equation}{section}
\begin{document}

\begin{center}

\end{center}

\centerline{\Large {\bf Embedding Dimensions }}
\bigskip
\centerline{\Large {\bf of Finite von Neumann Algebras}}

$\bigskip$

\centerline{  Junhao Shen}

\bigskip

\centerline{Mathematics Department, University of New Hampshire,
Durham, NH, 03824}

\bigskip

\centerline{ jog2@cisunix.unh.edu}

$\bigskip$

\noindent\textbf{Abstract: }  We introduce ``embedding dimensions"
of a family of generators of a finite von Neumann algebra  when the
von Neumann algebra can be faithfully embedded into the ultrapower
of the hyperfinite II$_1$ factor. These embedding dimensions are von
Neumann algebra invariants, i.e., do  not depend on the choices of
the generators. We also find values of these invariants for some
specific von Neumann algebras.

$\bigskip$

\section{Introduction}

 Let $H$ be a complex Hilbert
space and $B(H)$ the algebra of all bounded linear operators on $H$.
 A von Neumann algebra is a *-subalgebra of $B(H)$ that is
closed in the weak operator topology on $B(H)$.

In a   series of remarkable papers, published between 1936 and 1943,
  Murray and von Neumann  described some basic
structures on von Neumann algebras. They separated the family of
 von Neumann algebras into three types, I, II, and
III, and constructed examples for each type.
 In fact, Murray and von Neumann
 provided two methods for constructing
   type II$_1$ von Neumann algebras.
  One is  obtained from the ``left regular representation of a discrete
  infinite group."
 The other is related to the action of such a group on
 a measure space (of finite measure) by measure preserving transformations.

The first construction proceeds as follows. Let $H$ be $l^2(G)$. We
assume that $G$ is countable so that $H$
 is separable. For each $g$ in $G$, let $L_g$ be
  translation of functions in $l^2(G)$ by $g^{-1}$. Then $g
 \rightarrow L_g$ is a faithful unitary representation of $G$ on
 $H$. Let $L(G)$ be the von Neumann algebra generated by $\{ L_g: g \in
 G\}$. When  each conjugacy class in $G$ (other than that of the
 identity $e$) is infinite, $L(G)$
 is a factor of type II$_1$. In this case,
 we say that $G$ is an infinite conjugacy class (i.c.c.) group.

 Specific examples of such II$_1$ factors result from choosing for
 $G$ any of the free groups $F_n$ on $n$ generators ($n \ge 2$), or
 the
 direct products of two free groups $F_m\times F_p$, ($m,p\ge 2$),
 or the permutation group $\Pi$ of the integers $\Bbb Z$ (consisting
 of those permutations that leave fixed all but a finite subset of
 $\Bbb Z$).  A factor is called
 ``hyperfinite"  if it is the
 ultraweak closure of the ascending union of a family of
 finite-dimensional self-adjoint subalgebras.
 A deep result of Murray and von Neumann shows that all such factors
 are isomorphic. Moreover, $L(\Pi )$ is that hyperfinite factor of
 type II$_1$. Murray and von Neumann (\cite{MN}) show that $L(\Pi)$ is not *-isomorphic
 to $L(F_n)$ ($n \ge 2)$.

The second construction is more complicated.
   Let $(X, \mu)$ be a non-atomic measure space of finite measure. $G$ (with unit $e$)
 is a countable (infinite) group of measure preserving transformations of
  $X$. Our Hilbert space is $L^2(X, \mu)$. Let $\mathcal A$ be the commutative von Neumann algebra
 $L^{\infty}(X, \mu)$. Now $G$ can be viewed as a group of automorphisms
 of the von Neumann algebra $\mathcal A$. Murray and von Neumann constructed
  a von Neumann algebra $R(\mathcal A,G)$ associated with the group $G$ and
 the commutative von Neumann algebra $\mathcal A$. If   $G$ acts  freely and egodically on $X$,
  the
 von Neumann algebra $R(\mathcal A,G)$ is a type II$_1$ factor. In addition,
 $\mathcal A$ is a maximal abelian subalgebra in $R(\mathcal A,G)$. The normalizers of
 $\mathcal A$
 (those unitary operators $U$ in $R(\mathcal A,G)$ such that $U\mathcal AU^* =\mathcal A$)
  generate the von Neumann algebra $R(\mathcal A,G)$. We call
   a maximal abelian subalgebra of a finite von Neumann algebra,
   whose normalizers generate the full von Neumann algebra,
   a
 Cartan subalgebra.

It is a long standing open problem whether
  every
II$_1$
 von Neumann algebra has Cartan subalgebras.
This question was answered by Voiculescu  negatively after he
introduced his remarkable theory of   free entropy (see \cite {V2}),
an analogue of classical entropy and Fisher information measure.
Associated with the free entropy, he defined a free entropy
dimension $ \delta_0$ which, in some sense, measures the
``noncommutative dimension" of a space. In \cite {V2} he showed that
    for any $n$ in $\Bbb N$, $ \delta_0(L(  F_n)) \ge  n .$
      Soon he showed in \cite {V3})   that if a von Neumann algebra $\mathcal N$ has
     a Cartan subalgebra, then     $ \delta_0(\mathcal N)  \le 1 .$ Thus free
     group factors $L(F_n)$ ($n\ge 2$) have no Cartan subalgebra.
       Later, Ge (\cite {Ge2}) showed that if a von Neumann algebra $\mathcal N$ is
      not prime, i.e., is a tensor product of two infinite-dimensional von Neumann algebras, then
         $ \delta_0(\mathcal N)  \le 1 .$
      In particular, $L(  F_n) \ncong L(F_m)\otimes L(F_p)$ for all
      $n,m,p\ge 2$.

      In \cite{HaSh}, we introduced upper free
      orbit dimension for finite von Neumann algebras, a concept
      closely related to Voiculescu's free entropy dimension.  By some easily obtained properties of
      upper free orbit dimension, we got
      very general results which imply most of the applications of
     Voiculescu's free entropy dimension on finite von Neumann algebras.

It is well-known that Voiculescu's free entropy dimension is closely
related to Connes' embedding problem which asks whether every
separable type II$_1$ factor will be faithfully embedded into the
ultrapower of the hyperfinite II$_1$ factor, $\mathcal R^\omega$. In
fact Voiculescu's free entropy dimension of a finite von Neumann
algebra can be viewed as a measurement of the number of ways to
embed this von Neumann algebra into $\mathcal R^\omega$. The upper
free orbit-dimension of a von Neumann algebra, introduced in \cite
{HaSh}, can be view as a measurement of  the number of ways to
embed, modulo conjugate actions by unitary elements of $\mathcal
R^\omega$,
  this von Neumann algebra into $\mathcal R^\omega$.

In more details, suppose $x_1,\ldots, x_n$  is a family of
generators of a finite von Neumann $\mathcal N$ which can be
faithfully trace-preserving embedded into $\mathcal M_k(\Bbb
C)^\omega$, where $\mathcal M_k(\Bbb C)^\omega$ is the ultraproduct
of $\{\mathcal M_k(\Bbb C)\}_{k=1}^\infty$ along the free filter
$\omega$. Based on the  philosophy in preceding paraghaph, we define
$\mathcal H_s^\omega(x_1,\ldots,x_n) $, $s$-embedding dimension of
$x_1,\ldots, x_n$ for $s\ge 0$, to be some measurement of the number
of ways to embed $\mathcal N$ into the ultrapower $\mathcal M_k(\Bbb
C)^\omega$. We show that $\mathcal H_s^\omega$ is a von Neumann
algebra invariant, i.e. does not depend on the choices of families
of generators. Then we carry out the computation of values of these
invariants for finite von Neumann algebras. For example, if
$\mathcal N$ is a hyperfinite von Neumann algebra, then $\mathcal
H_0^\omega(\mathcal N) =0$. If $\mathcal N$ is the free group factor
on $n$ generators, then $\mathcal H_1^\omega(\mathcal N)=\infty$. On
the other hand, $\mathcal H_1^\omega(\mathcal N)=0$ if $\mathcal N$
is a type II$_1$ factor with Cartan subalgebras, a nonprime type
II$_1$ factor , or some type II$_1$ factor with property T.
Therefore, this invariant does give us   information on the
classification  of type II$_1$ factors.

  Because these ``embedding dimensions" may have
potential applications on Connes' embedding problems, it is
worthwhile to study them in more details. Having its motivations
from Voiculescu's free entropy dimension and from free orbit
dimension of \cite{HaSh}, we develop the theory of embedding
dimensions from its own interest and keep the paper as
self-contained as possible. Anther motivation of the paper comes
from the attempt to further classify II$_1$ von Neumann algebras
whose Voiculescu's free entropy dimensions are equal to $1$,
especially from the question whether the tensor products of free
group factors have Cartan subalgebras. We can view $\mathcal
H_s^\omega(\mathcal N)$, $s$-embedding-dimension of finite von
Neumann algebra $\mathcal N$, as an analogue of the classical
fractal dimension  in the subject of finite von Neumann algebras. It
is not hard to see that $\mathcal H_s^\omega(\mathcal N)$ is a
decreasing function of $s\ge 0$. We know that $\mathcal
H_1^\omega(\mathcal N)=0$ for many type II$_1$ factors (for example,
type II$_1$ factors with Cartan subalgebras). Hopefully  $\mathcal
H_s^\omega(\mathcal N)$ becomes non zero when $s$ is small enough.
In this direction we wish that these ``fractal" dimensions can
provide us with new tools to further classify type II$_1$ factors
whose Voiculescu's free entropy dimensions are equal to $1$.

The organization of the paper is as follows. In section 2, we
introduce some notations and give the definitions of embedding
dimensions. We show that these embedding dimensions are von Neumann
algebra invariants in section 3. The embedding dimensions of abelian
von Neumann algebras and free group factors are obtained in section
4. Some technical  lemmas on covering numbers and unitary orbit
covering numbers are proved in section 5. The computation of values
of embedding dimension for some specific type II$_1$ factors is
carried out in section 6.

\section{Some notations and definitions}

\subsection{Covering numbers}  Let $\mathcal{M}_{k}(\mathbb{C})$ be the $k\times k$
full matrix algebra with entries in $\mathbb{C}$, and $\tau_{k}$ be
the normalized trace on $\mathcal{M}_{k}(\mathbb{C})$, i.e.,
$\tau_{k}=\frac{1}{k}Tr$, where $Tr$ is the usual trace on
$\mathcal{M}_{k}(\mathbb{C})$. Let $\mathcal{U}(k)$ denote the group
of all unitary matrices in $\mathcal{M}_{k}(\mathbb{C})$. Let
$\mathcal{M}_{k}(\mathbb{C})^{n}$ denote the direct sum of $n$
copies of $\mathcal{M}_{k}(\mathbb{C})$. Let $\Vert\cdot\Vert_{2}$
denote the trace norm induced by $\tau_{k}$ on
$\mathcal{M}_{k}(\mathbb{C})^{n}$, i.e.,
\[
\Vert(A_{1},\ldots,A_{n})\Vert_{2}^{2}=\tau_{k}(A_{1}^{\ast}A_{1})+\ldots
+\tau_{k}(A_{n}^{\ast}A_{n})
\]
for all $(A_{1},\ldots,A_{n})$ in $\mathcal{M}_{k}(\mathbb{C})^{n}$.
\begin{definition}
For every $\delta>0$, we define the $\delta$-ball $Ball(B_{1},\ldots
,B_{n};\delta)$ centered at $(B_{1},\ldots,B_{n})$ in
$\mathcal{M}_{k}(\mathbb{C})^{n}$ to be the subset of
$\mathcal{M}_{k}(\mathbb{C})^{n}$   consisting of all
$(A_{1},\ldots,A_{n})$ in $\mathcal{M}_{k}(\mathbb{C})^{n}$ such
that
$\Vert(A_{1},\ldots,A_{n})-(B_{1},\ldots,B_{n})\Vert_{2}<\delta.$
\end{definition}
\begin{definition}
For every $\delta>0$, we define the $\delta$-orbit-ball $\mathcal{U}%
(B_{1},\ldots,B_{n};\delta)$ centered at $(B_{1},\ldots,B_{n})$ in $\mathcal{M}%
_{k}(\mathbb{C})^{n}$ to be the subset of
$\mathcal{M}_{k}(\mathbb{C})^{n}$ consisting of all
$(A_{1},\ldots,A_{n})$ in $\mathcal{M}_{k}(\mathbb{C})^{n}$ such
that there exists some unitary matrix $W$ in $\mathcal{U}(k)$
satisfying
\[
\Vert(A_{1},\ldots,A_{n})-(WB_{1}W^{\ast},\ldots,WB_{n}W^{\ast})\Vert
_{2}<\delta.
\]
\end{definition}

\begin{definition}For every $R>0$, denote by $(\mathcal M_k(\Bbb C)^n)_R$ the
subset of $\mathcal M_k(\Bbb C)^n$ consisting of all these
$(A_1,\ldots,A_n)$ satisfying $\max_{1\le j\le n}\| A_j\|\le R$.

\end{definition}

\begin{definition} Let $\Gamma$ be a subset of  $\mathcal{M}%
_{k}(\mathbb{C})^{n}$. (i) For $\delta>0$, we define the {\em
$\delta$-covering number} $\nu_2(\Gamma,\delta)$ to be the minimal
number of $\delta$-balls that cover $\Gamma$ with the centers of
these $\delta$-balls in $\Gamma$. (ii) Define the {\em
$\delta$-orbit covering number} $\nu(\Gamma,\delta)$ to be the
minimal number of $\delta$-orbit-balls that cover $\Gamma$ with the
centers of these $\delta$-orbit-balls in $\Gamma$.

\end{definition}

\subsection{Embedding dimensions}
 Let ${\mathcal M_k(\Bbb C)}$ be the $k\times k$ full matrix algebra with complex entries and
$\tau_k$ be the normalized trace on $\mathcal M_k(\Bbb C)$. If
$\omega$ is a free filter in $\beta(\Bbb N)\setminus \Bbb N $  then
denote by $\mathcal M_k(\Bbb C)^{\omega}$ the quotient of the von
Neumann algebra $l^{\infty}(\Bbb N, \prod_{k=1}^\infty\mathcal
M_k(\Bbb C))$ by the $0$-ideal of the trace $\tau_{\omega}$, where
$\tau_\omega$ is defined by  $ \tau_{\omega}
((X_k)_k)=\lim_{k\rightarrow \omega}\tau_k(X_k)$ for every
$x=(X_k)_k$ in $ \mathcal M_k(\Bbb C )^{\omega}$. We also define the
Hilbert norm $\| \ \|_{2,\omega}$ on $ \mathcal M_k(\Bbb C )^\omega$
by
$$\|x||_{2,\omega}=\|(X_k)_k\|_{2,\omega}=\tau_\omega([(X_k)_k]^*[(X_k)_k])^{1/2}=\lim_{k\rightarrow
\omega} \|X_k\|_2$$ for all $x=(X_k)_k$ in $ \mathcal M_k(\Bbb C
)^\omega$. Then $ \mathcal M_k(\Bbb C )^{\omega}$ is a type II$_1$
factor; and $\tau_{\omega}$ is a tracial trace on $\mathcal M_k(\Bbb
C )^{\omega}$.

Let $\mathcal N$ be a finitely generated von Neumann algebra with a
tracial state $\tau$. Assume that $\mathcal N$ can be faithfully
embedded into $\mathcal M_k(\Bbb C)^{\omega}$.

\begin{definition} We have the following definitions.

(i) Define $\Theta(\mathcal N, \mathcal M_k(\Bbb C)^{\omega})$ as a
subset of $Hom(\mathcal N, \mathcal M_k(\Bbb C)^{\omega})$
consisting of all faithful trace-preserving embedding $\theta$ from
$(\mathcal N,\tau)$ into $(\mathcal M_k(\Bbb
C)^{\omega},\tau_{\omega}).$

(ii) Define $\Xi$ to be the set consisting of all sequences, $\xi=
\{k_m\}_{m=1}^\infty$,  of positive integers  such that
$\lim_{m\rightarrow \infty} k_m=\infty.$

(iii) We will also introduce   ``Voiculescu's topological structure"
on the space
$$\mathfrak T= \left ( \prod_{k=1}^\infty \mathcal M_k(\Bbb C)\right)^n=\prod_{k=1}^\infty \mathcal M_k(\Bbb C)\times \cdots\times
\prod_{k=1}^\infty \mathcal M_k(\Bbb C).$$ First we introduce the
neighborhood in $\mathfrak T$, indexed by   elements $y_1,\ldots,
y_n$ in $\mathcal M_k(\Bbb C)^{\omega}$,   $R>0$ and
$\xi=\{k_m\}_{m=1}^\infty$ in $\Xi$ as follows. Define the
neighborhood
$$\mathfrak N_{R,\xi}(y_1,\ldots,y_n)$$ as a subset of $\prod_{k=1}^\infty \mathcal M_k(\Bbb
C)\times \cdots\times \prod_{k=1}^\infty \mathcal M_k(\Bbb C)$
consisting of all   $$ ((Y_{1,k})_{k=1}^\infty,\ldots,
(Y_{n,k})_{k=1}^\infty)$$   in $\prod_{k=1}^\infty \mathcal M_k(\Bbb
C)\times \cdots\times \prod_{k=1}^\infty \mathcal M_k(\Bbb C)$ such
that $\|Y_{i,k}\|< R$ for all $1\le i\le n, k\ge 1$ and
$$ |\tau_{k}(Y_{j_1}^{\epsilon_1}\ldots Y_{j_p}^{\epsilon_p})
-\tau_{\omega}(y_{j_1}^{\epsilon_1}\ldots
y_{j_p}^{\epsilon_p})|<\frac 1 m,$$ for all $k\ge k_m$,  $1\le p\le
m$, $1\le j_1,\ldots, j_p\le n$ and $\epsilon_i\in \{*,1\}$ for
$1\le i\le p.$

(iv) For each $k\ge 1$, let $\Bbb P_k$ be the projection from
$(\prod_{k=1}^\infty\mathcal M_k(\Bbb C))^n$ onto $(\mathcal
M_k(\Bbb C))^n.$
\end{definition}

\begin{remark} The introductions of the set $\Xi$ and the ``topology" are very necessary.
Because the first finitely many terms of any representative of an
element in $\mathcal M_k(\Bbb C)^\omega$ can be chosen  arbitrarily,
we use this $\Xi$ to exclude this ``arbitrary" phenomena. Each $\xi$
in $\Xi$ plays the role of radius.  And it is not hard to see that
$$ \bigcup_{ \{y_1,\ldots, y_n\}\subset \mathcal M_k(\Bbb C)^{\omega}}\left
(\bigcup_{R>0, \ \xi\in\Xi}\ \mathfrak N_{R,\xi}(y_1,\ldots,y_n)
\right)= l^\infty(\Bbb N, \prod_{k=1}^\infty \mathcal M_k(\Bbb
C)\times \cdots\times \prod_{k=1}^\infty \mathcal M_k(\Bbb C)).
$$

\end{remark}

\begin{definition}  Suppose that $x_1,\ldots,x_n$ is a family of
elements in $\mathcal N$. Let $R$ be a positive number and
$\xi=\{k_m\}_{m=1}^\infty$ be in $\Xi$. We define the neighborhood
in $\mathfrak T$, indexed by $x_1,\ldots,x_n$ in $\mathcal N$, $R>0$
and $\xi\in\Xi$, by
$$ \mathfrak N_{R,\xi}(x_1,\ldots, x_n) =\bigcup_{\theta\in\Theta(\mathcal N,\mathcal M_k(\Bbb C)^\omega)} \ \mathfrak N_{R,\xi}(\theta(x_1),\ldots,\theta(x_n)).   $$

{\em Voiculescu's embedding dimension}  of $x_1,\ldots,x_n$,
$\delta_0^\omega (x_1,\ldots,x_n)$, is defined by
\[
\begin{aligned}
 \delta_0^\omega (x_1,\ldots, x_n;R,\xi,\delta)& =\  \lim_{k\rightarrow \omega} \frac
{\log(\nu_2(\Bbb P_k(\mathfrak N_{R,\xi}( x_1,\ldots,
 x_n)),\delta))}{-k^{2}\log\delta}\\
\delta_0^\omega (x_1,\ldots, x_n)&= \limsup_{\delta\rightarrow 0}\
\sup_{R>0,\xi\in \Xi}  \delta_0^\omega (x_1,\ldots,
x_n;R,\xi,\delta)
\end{aligned}
\]
 For   $s\ge 0$, we define {\em $s$-embedding dimension}  of $x_1,\ldots,x_n$, $\mathcal H_s^\omega  (x_1,\ldots,
 x_n)$,  by\[
\begin{aligned}
\mathcal H_s^\omega  (x_1, \ldots,
 x_n;R,\xi,\delta) &=  \lim_{k\rightarrow \omega} \frac
{\log(\nu(\Bbb P_k(\mathfrak N_{R,\xi}( x_1,\ldots,
 x_n)), \delta))}{k^{2s}}\\
\mathcal H_s^\omega  (x_1, \ldots,
 x_n) &=
\limsup_{\delta\rightarrow 0}\ \sup_{R>0,\xi\in \Xi}  \mathcal
H_s^\omega  (x_1, \ldots,
 x_n;R,\xi,\delta)
\end{aligned}
\]
Generally, given any function $f(s, \cdot)$ where $s$ is a parameter
or a family of parameters, we define $f(s, \cdot)$-dimension of
$x_1,\ldots,x_n$, $\mathcal H_{f(s,\cdot)}^\omega  (x_1,\ldots,
 x_n)$, by
 \[
\begin{aligned}\mathcal H_{f(s,\cdot)}^\omega  (x_1, \ldots,
 x_n;R,\xi,\delta) &=  \lim_{k\rightarrow \omega} \frac
{\log(\nu(\Bbb P_k(\mathfrak N_{R,\xi}( x_1,\ldots,
 x_n)), \delta))}{f(s,k)}\\
\mathcal H_{f(s,\cdot)}^\omega  (x_1, \ldots,
 x_n) &=
\limsup_{\delta\rightarrow 0}\ \sup_{R>0,\xi\in \Xi}  \mathcal
H_{f(s,\cdot)}^\omega  (x_1, \ldots,
 x_n;R,\xi,\delta)
\end{aligned}
\]
\end{definition}

\begin{remark}
It follows from the definition that $\mathcal
H_s^\omega(x_1,\ldots,x_n;R, \xi,\delta)$ is a decreasing function
of $\delta>0$ and $\mathcal H_s^\omega(x_1,\ldots,x_n)$ is a
decreasing function of $s\ge 0$.
\end{remark}

\begin{remark}
When $s>1$, we have $\mathcal H_s^\omega(x_1,\ldots,x_n)=0$  for all
$x_1,\ldots,x_n$ in $\mathcal N$.
\end{remark}

\begin{remark}  Voiculescu's embedding dimension
could be viewed as a measurement of the number of the ways to  embed
$\mathcal N$ into $\mathcal M_k(\Bbb C)^\omega$. For every embedding
$\theta$ from $\mathcal N$ into $\mathcal M_k(\Bbb C)^\omega$, we
know that $u\theta(\cdot)u^*$ is also an embedding from $\mathcal N$
into $\mathcal M_k(\Bbb C)^\omega$, where $u$ is a unitary element
of $\mathcal M_k(\Bbb C)^\omega$. Therefore the
``$s$-embedding-dimension", or $\mathcal H_s$, could be viewed as a
measurement of the number of the ways to embed $\mathcal N$, modulo
conjugate actions by  unitary elements of $\mathcal M_k(\Bbb
C)^\omega$, into $\mathcal M_k(\Bbb C)^\omega$.
\end{remark}

We should also define the embedding dimensions of a family of
elements $x_{1},\ldots,x_{n}$ in the presence of another family of
elements $y_{1},\ldots,y_{p}$ of $\mathcal N$.

\begin{definition} Suppose \
$x_1,\ldots,x_n,y_1,\ldots,y_p$ \ are the elements of $\mathcal N$.
\ Let $R>0$, $\xi$ be in $\Xi$ and $\theta$ be in $\Theta(\mathcal
N,\mathcal M_k(\Bbb C)^\omega)$. Let
$$\mathfrak N_{R,\xi}(\theta(x_{1}),\ldots,\theta(x_{n}):\theta(y_{1})%
,\ldots,\theta(y_{p}))$$ be the image of the projection of
$\mathfrak N_{R,\xi}(\theta(x_{1}),\ldots,\theta(x_{n}),\theta(y_{1})%
,\ldots,\theta(y_{p}))$ onto the first $n$ components, i.e.,
\[
((A_{1,k})_k,\ldots,(A_{n,k})_k)\in\mathfrak
N_{R,\xi}(\theta(x_{1}),\ldots,\theta(x_{n}):\theta(y_{1})%
,\ldots,\theta(y_{p}))
\]
if there are elements $B_{1,k},\ldots,B_{p,k}$ in
$\mathcal{M}_{k}(\mathbb{C})$ such that
\[
((A_{1})_k,\ldots,(A_{n,k})_k,(B_{1,k})_k,\ldots,(B_{p,k})_k)\in\mathfrak N_{R,\xi}(\theta(x_{1}),\ldots,\theta(x_{n}), \theta(y_{1})%
,\ldots,\theta(y_{p})).
\] Then we define,
$$\begin{aligned}
\mathfrak N_{R,\xi}(x_1,\ldots, x_n: & \ y_1,\ldots,y_p)
\\ &=\bigcup_{\theta\in\Theta(\mathcal N,\mathcal M_k(\Bbb C)^\omega)}
\ \mathfrak N_{R,\xi}(\theta(x_1),\ldots,\theta(x_n):
\theta(y_1),\ldots, \theta(y_p)),\\
 \delta_0^\omega  (x_1,\ldots, x_n:&\ y_1,\ldots,y_p; R,\xi,\delta)\\ & =   \lim_{k\rightarrow \omega} \frac
{\log(\nu_2(\Bbb P_k(\mathfrak N_{R,\xi}( x_1,\ldots,
 x_n:y_1,\ldots,y_p)),\delta))}{-k^{2}\log\delta}\\
 \delta_0^\omega  (x_1,\ldots, x_n:&\ y_1,\ldots,y_p)\\ & =\limsup_{\delta\rightarrow 0}\ \sup_{R>0,\xi\in\Xi}   \delta_0^\omega  (x_1,\ldots, x_n: \ y_1,\ldots,y_p; R,\xi,\delta)
\end{aligned}$$
 And $$\mathcal H_s^\omega (x_1,\ldots,x_n:y_1,\ldots,y_p) , \quad \mathcal
H_{f(s,\cdot)}^\omega (x_1,\ldots,x_n:y_1,\ldots,y_p)$$ are defined
similarly.
\end{definition}

\section{$\mathcal H_s^\omega$ is a von Neumann algebra invariant}
In this section, we are going to show that $\mathcal H_s^\omega$ is
a von Neumann algebra invariant, i.e. it does not depend on the
choices of the generators. First, we have the following lemma which
follows directly from the definition of embedding dimensions.

\begin{lemma}Suppose $\mathcal N$ is a finitely generated von
Neumann algebra with a tracial state $\tau$ and $\mathcal N$ can be
faithfully trace-preserving embedded into $\mathcal M_k(\Bbb
C)^\omega$. Let $x_{1},\ldots,x_{n},y_{1},\ldots,y_{p}$ be elements
in a von Neumann algebra $\mathcal{N}$. If $x_{1},\ldots,x_{n}$
generate $\mathcal N$ as a von Neumann algebra,  then, for every
$s\ge 0$,
\[
\mathcal{H}_s^\omega(x_{1},\ldots,x_{n})=\mathcal{H}_s^\omega(x_{1},\ldots,x_{n}%
:y_{1},\ldots,y_{p})\mathfrak{.}%
\]

\end{lemma}
\begin{proof}
Note that $x_1,\ldots,x_n$ generate $\mathcal N$ and
$y_1,\ldots,y_p$ are contained in $\mathcal N$. Thus  $y_1,\ldots,
y_p$ can be approximated by the polynomials of $x_1,\ldots, x_n$ in
$\| \cdot \|_2$-norm. Let $R>\max\{\|x_i\|, \|y_j\|,1\le i\le n,
1\le j\le p\}$ and $\xi=\{k_m\}_{m=1}^\infty$ be in $ \Xi$.  For
each $m\ge 0$, there is $m'\ge 0$ satisfying: for every $$
\{X_1,\ldots, X_n,Y_1,\ldots, Y_p\}\subset \mathcal (M_k(\Bbb
C))_R,$$ if
$$ |\tau_{k}(X_{j_1}^{\epsilon_1}\ldots X_{j_q}^{\epsilon_q})
-\tau (x_{j_1}^{\epsilon_1}\ldots x_{j_q}^{\epsilon_q})|<\frac 1
{m'},$$ for all $\{X_{j_i}\}_{i=1}^q\subset \{X_1,\ldots,X_n\}$,
$\{\epsilon_i\}_{i=1}^q\subset \{*,1\}$  and $1\le q\le m'$, then
$$ |\tau_{k}(Z_{j_1}^{\epsilon_1}\ldots Z_{j_q}^{\epsilon_q})
-\tau (z_{j_1}^{\epsilon_1}\ldots x_{j_q}^{\epsilon_q})|<\frac 1 {m
},$$ for all   $\{Z_{j_i}\}_{i=1}^q\subset \{X_1,\ldots, X_n,
Y_1,\ldots,Y_p\}$, $\{z_{j_i}\}_{i=1}^q\subset \{x_1,\ldots, x_n,
y_1,\ldots,y_p\}$,
 $\{\epsilon_i\}_{i=1}^q\subset \{*,1\}$ and  $1\le q\le m$.

Let $\tilde \xi=\{\tilde k_{m'}\}_{m'=1}^\infty$ in $\Xi$ such that
$\tilde k_{m'}=k_m$. Then
$$\begin{aligned}
\mathfrak N_{R,\tilde\xi}(x_1,\ldots,x_n) \subset\mathfrak
N_{R,\xi}(x_1,\ldots,x_n:y_1,\ldots,y_p) \subset \mathfrak N_{R,
\xi}(x_1,\ldots,x_n),\end{aligned}
$$ for all $k\ge 1$. The rest follows from the definitions.
\end{proof}

Now we are ready to show the main result in this section.

\begin{theorem}
Suppose $\mathcal{N}$ is a von Neumann algebra with a tracial state
$\tau$ and $\mathcal N$ can be faithfully trace-preserving embedded
into $\mathcal M_k(\Bbb C)^\omega$. Suppose $\{x_{1},\ldots,x_{n}\}$
and $\{y_1,\ldots, y_p\}$ are two families of generators of
$\mathcal N$. Then, for all $s\ge 0$,
\[
\mathcal H_s^\omega(x_1,\ldots,x_n)=\mathcal
H_s^\omega(y_1,\ldots,y_p).
\]

\end{theorem}

\begin{proof}
\textbf{\ } We  need only to show that
\[
\mathcal H_s^\omega(x_1,\ldots,x_n)\ge \mathcal
H_s^\omega(y_1,\ldots,y_p).
\]

Since $x_{1},\ldots,x_{n}$ are elements in $\mathcal{N}$ that
generate $\mathcal{N}$ as a von Neumann algebra, for every
$0<\delta<1$,
there exists a family of noncommutative polynomials $\psi_{i}(x_{1}%
,\ldots,x_{n})$, $1\leq i\leq p$, such that
\[
\sum_{i=1}^{p}\Vert
y_{i}-\psi_{i}(x_{1},\ldots,x_{n})\Vert_{2}^{2}<\left(
\frac{\delta}{4}\right)  ^{2}.
\]
Therefore, for every $\theta$ in $\Theta(\mathcal N,\mathcal
M_k(\Bbb C)^\omega)$, we have
\[
\sum_{i=1}^{p}\Vert
\theta(y_{i})-\psi_{i}(\theta(x_{1}),\ldots,\theta(x_{n}))\Vert_{2}^{2}<\left(
\frac{\delta}{4}\right)  ^{2}.
\]
 For such a family of polynomials
$\psi_{1},\ldots,\psi_{p}$, and every $R>0$
there always exists a constant $D\geq1$, depending only on $R,\psi_{1}%
,\ldots,\psi_{n}$, such that
\[
\left (
\sum_{i=1}^{p}\Vert\psi_{i}(A_{1},\ldots,A_{n})-\psi_{i}(B_{1},\ldots
,B_{n})\Vert_{2}^{2}\right )^{1/2}\leq D
\Vert(A_{1},\ldots,A_{n})-(B_{1},\ldots ,B_{n})\Vert_{2},
\]
for all $(A_{1},\ldots,A_{n}),(B_{1},\ldots,B_{n})$ in $\mathcal{M}%
_{k}(\mathbb{C})^{ {n}}$, all $k\in\mathbb{N}$, satisfying $\Vert
A_{j}\Vert, \Vert B_{j}\Vert\leq R,$ for $1\leq j\leq n.$

For $R> 1, \xi\in\Xi, $  by the definition of $\mathfrak
N_{R,\xi}(\theta(y_1),\ldots,
\theta(y_p),\theta(x_1),\ldots,\theta(x_n))$, there is some $m_0$
(only depending on $\delta$, not on $\theta$) such that, when $k\ge
k_{m_0}$, every
$$(\tilde H_{1},\ldots,\tilde  H_{p},\tilde A_{1},\ldots,\tilde A_{n}) \quad \in \quad
 \Bbb P_k(\mathfrak
N_{R,\xi}(\theta(y_1),\ldots,
\theta(y_p),\theta(x_1),\ldots,\theta(x_n)))$$ satisfies
$$\begin{aligned}
&\left |  \left (\sum_{i=1}^{p}\Vert \tilde H_{i}- \psi_{i}(\tilde
A_{1},\ldots,\tilde A_{n})\Vert_{2}^{2}\right)^{1/2} \right .
\\ &\qquad \qquad \left .-  \left (\sum_{i=1}^{p}\Vert
\theta(y_{i})-\psi_{i}(\theta(x_{1}),\ldots,\theta(x_{n})\Vert_{2}^{2}\right
)^{1/2} \right | &\leq \frac{\delta}{4}\end{aligned}
$$ Hence,
\[
\left (\sum_{i=1}^{p}\Vert \tilde H_{i}-\psi_{i}(\tilde
A_{1},\ldots,\tilde A_{n})\Vert_{2}^{2}\right )^{1/2}\leq
\frac{\delta}{2} .
\] It follows that, when $k\ge k_{m_0}$, every
$$(H_{1},\ldots,H_{p},A_{1},\ldots,A_{n}) \quad \in \quad
 \Bbb P_k(\mathfrak
N_{R,\xi}(y_1,\ldots,  y_p,x_1,\ldots,x_n))$$ satisfies
\[
\left (\sum_{i=1}^{p}\Vert
H_{i}-\psi_{i}(A_{1},\ldots,A_{n})\Vert_{2}^{2}\right )^{1/2}\leq
\frac{\delta}{2} .
\]
 It is
obvious that such an $(A_{1},\ldots,A_{n})$ is also in $\Bbb
P_k(\mathfrak N_{R,\xi}(x_1,\ldots,x_n: y_1,\ldots, y_p))$, which is
contained in $\Bbb P_k(\mathfrak N_{R,\xi}(x_1,\ldots,x_n))$. On the
other hand, by the definition of the orbit covering number, we know
there exists a set $\{ \mathcal{U}(B_{1}^{\lambda
},\ldots,B_{n}^{\lambda};\frac{\delta}{4D})
\}_{\lambda\in\Lambda_{k}}$ of $\frac{\delta}{4D}$-orbit-balls that
cover $\Bbb P_k(\mathfrak N_{R,\xi}(x_1,\ldots,x_n))$ with the
cardinality of $ \Lambda_{k} $ satisfying \
$|\Lambda_{k}|=\nu(\Bbb P_k(\mathfrak N_{R,\xi}(x_1,\ldots,x_n)),\frac{\delta}{4D}).$ \ Thus \ for \ such \ $(A_{1}%
,\ldots,A_{n})$ in $\Bbb P_k(\mathfrak N_{R,\xi}(x_1,\ldots,x_n))$,
there exists some $\lambda\in\Lambda_{k}$, such
that
$$
  (A_1,\ldots, A_n) \in \mathcal{U}(B_{1}^{\lambda
},\ldots,B_{n}^{\lambda};\frac{\delta}{4D}),
$$ i.e., there is some
$W\in\mathcal{U}(k)$ such that
\[
\Vert(A_{1},\ldots,A_{n})-(WB_{1}^{\lambda}W^{\ast},\ldots,WB_{n}^{\lambda
}W^{\ast})\Vert_{2}\leq\frac{\delta}{4D}.
\]
It follows that
\[
\sum_{i=1}^{p}\Vert
H_{i}-W\psi_{i}(B_{1}^{\lambda},\ldots,B_{n}^{\lambda
})W^{\ast}\Vert_{2}^{2}=\sum_{i=1}^{p}\Vert
H_{i}-\psi_{i}(WB_{1}^{\lambda
}W^{\ast},\ldots,WB_{n}^{\lambda}W^{\ast})\Vert_{2}^{2}\leq
 {\delta}   ^{2},
\]
for some $\lambda\in\Lambda_{k}$ and $W\in\mathcal U(k),$ i.e.,
\[
(H_{1},\ldots,H_{p})\in\mathcal{U}(\psi_{1}(B_{1}^{\lambda},\ldots
,B_{n}^{\lambda}),\ldots,\psi_{p}(B_{1}^{\lambda},\ldots,B_{n}^{\lambda
}); \delta).
\]
Hence, by the definition of   embedding dimension, we get
\[
\begin{aligned}
0&\le \mathcal H_s^\omega(y_1,\ldots,y_p:
x_1,\ldots,x_n;R,\xi,2\delta)\le
 \lim_{k\rightarrow \omega} \frac
{\log(|\Lambda_k|)}{ k^{2s}}\\
&=
  \lim_{k\rightarrow \omega} \frac
{\log(\nu(\Bbb P(\mathfrak N_{R,\xi}(x_1,\ldots, x_n )),\frac
\delta{4D}))}{ k^{2s} }\\
& \le  \mathcal H_s^\omega(x_{1},\ldots,x_{n}),
\end{aligned}
\]
since $$ \mathcal H_s^\omega(x_{1},\ldots,x_{n}; R,\xi,\delta_1)
=\sup_{\delta_1>0}\ \sup_{R>0,\xi\in\Xi} \ \lim_{k\rightarrow
\omega} \frac {\log(\nu(\Bbb P(\mathfrak N_{R,\xi}(x_1,\ldots, x_n
)),\delta_1))}{ k^{2s} }.
$$Therefore $\mathcal H_s^\omega
(y_{1},\ldots,y_{p}:x_{1},\ldots,x_{n})\le \mathcal
H_s^\omega(x_{1},\ldots,x_{n})$. Now it follows from Lemma 1 that
\[
\mathcal H_s^\omega (y_{1},\ldots,y_{p})=\mathcal H_s^\omega
(y_{1},\ldots ,y_{p}:x_{1},\ldots,x_{n}).
\] Hence $\mathcal H_s^\omega (y_{1},\ldots,y_{p})\le\mathcal H_s^\omega(x_{1},\ldots,x_{n})$, which completes the proof.
\end{proof}

Because of the preceding theorem, the following definition is
well-defined.
\begin{definition}
Suppose $\mathcal N$ is a finitely generated von Neumann algebra
with a tracial state $\tau$ and $\mathcal N$ can be faithfully
trace-preserving embedded into $\mathcal M_k(\Bbb C)^\omega$ for
some free  filter $\omega$ in $\beta(\Bbb N)\setminus \Bbb N$. Then,
we define the $s$-embedding-dimension of $\mathcal N$ by
$$
\mathcal H_s^\omega(\mathcal N)=\mathcal H_s^\omega(x_1,\ldots,x_n),
$$ for any family of generator $x_1,\ldots,x_n$ of $\mathcal N$.
\end{definition}

\begin{remark}
 It is trivial to see that $\mathcal H_s^\omega(\mathcal N)$ is a
 decreasing function of $s\ge 0$ for every finite von Neumann algebra $\mathcal
 N$. And $\mathcal H_s^\omega(\mathcal N)=0$ for all $s>1$.
\end{remark}
\begin{remark}
We did not prove the Voiculescu's embedding dimension is a von
Neumann algebra invariant.
\end{remark}

\section{Embedding dimensions of abelian von Neumann algebras and
free group factors}

In this section, we are going to compute values of the embedding
dimensions of abelian von Neumann algebras and free group factors.

The following lemma was first proved by Voiculescu in \cite{V2}. A
simplified proof can be found in \cite{DH}. For the sake of
completeness, we also sketch its proof here.

\begin{lemma}
Let $x$ be a self-adjoint element in a von Neumann algebra $\mathcal
N$ with a tracial state $\tau$. Let $R>\| x \|$. For every
$\delta>0$, there is some positive integer $m$ such that, for all
$k\ge 1$, if $A, B$ are two self-adjoint matrices in $\mathcal
M_k(\Bbb C)$ satisfying $\|A\|\le R, \|B\|\le R$ and
$$
\begin{aligned}
  | \tau_k(A^p) -\tau(x^p)| < \frac 1 m; \qquad | \tau_k(B^p) -\tau(x^p)| < \frac 1 m,
\end{aligned}
$$ for all $1\le p\le m$, then there is some unitary matrix $U$ in
$\mathcal U(k)$ such that
$$ \| UAU^*-B\|_2\le \delta.  $$
\end{lemma}

\begin{proof}
Suppose on the contrary that the following holds: there is some
$\delta_0>0$ such that
 for every $m\ge 1$, there is some $k_m\ge 1$ and
some self-adjoint matrices $A_m, B_m$ in $\mathcal M_{k_m}(\Bbb C)$
satisfying $\|A_m\|\le R, \|B_m\|\le R$,
$$
\begin{aligned}
  | \tau_{k_m}(A_m^p) -\tau(x^p)| < \frac 1 m; \qquad | \tau_{k_m}(B_m^p) -\tau(x^p)| < \frac 1 m,
\end{aligned}
$$ for all $1\le p\le m$, and  $ \| UA_mU^*-B_m\|_2> \delta_0 $ for all  unitary matrix $U$ in
$\mathcal U(k_m)$.

 Let $\omega$ be a free filter in $\beta(\Bbb N)\setminus\Bbb N$. Denote by $\mathcal
M_{k_m}(\Bbb C)^{\omega}$ the ultrapower of $\{\mathcal M_{k_m}(\Bbb
C)\}_{m=1}^\infty$ along the  filter $\omega$.  It is not hard to
see $[(A_m)_m], [(B_m)_m]$ are two self-adjoint elements in
$\mathcal M_{k_m}(\Bbb C)^{\omega}$ that have the same distribution.
By Lemma 7.1 of \cite {Po1}, there is some $u$ in $\mathcal
M_{k_m}(\Bbb C)^{\omega}$, such that $u[(A_m)_m] u^*= [(B_m)_m]$.
Let $(U_m)_m$ be a representative of $u$ in $\mathcal M_{k_m}(\Bbb
C)^{\omega}$. We can assume that each $U_m$ is a unitary matrix in
$\mathcal M_{k_m}(\Bbb C)$. Hence $\lim_{m\rightarrow
\omega}\|U_mA_mU_m^*-B_m\|_2=0$, which contradicts with the
assumption that $ \| UA_mU^*-B_m\|_2> \delta_0 $ for all unitary
matrix $U$ in \ $\mathcal U(k_m)$. Therefore, the statement of the
lemma is true.

\end{proof}

\begin{remark}
The proof of the preceding lemma  shows the same statement also
holds for a unitary element $x$ in $\mathcal N$ (considering
$*$-distribution of a unitary element instead of  distribution of
 a self-adjoint element). In fact, a stronger result was obtained in
\cite{DH} in the case when $x_1,\ldots,x_n$ generate a hypefinite
von Neumann algebra.
\end{remark}
\begin{theorem}
Suppose $\mathcal A$ is an abelian von Neumann algebra with a
tracial state $\tau$. Then $\mathcal H_0^\omega(\mathcal A)=0$.

\end{theorem}

\begin{proof}
By \cite{MN}, we can assume that the abelian von Neumann algebra
$\mathcal A$ is generated by a self-adjoint element  $x$. It is
well-known that every abelian von Neumann algebra with a tracial
state can be faithfully trace-preserving embedded into the
ultrapower $\mathcal M_k(\Bbb C)^\omega$. Let $\delta>0$, $R>0$ and
$\xi=\{k_m\}_{m=1}^\infty$ be in $\Xi$. For every
$$
A_k,  B_k \in \Bbb P_k(\mathfrak N_{R,\xi}(x)),
$$ there are some $\theta_1,\theta_2 $ in $\Theta (\mathcal A, \mathcal M_k(\Bbb
C)^\omega)$ such that
$$
A_k= \Bbb P_k(\mathfrak N_{R,\xi}(\theta_1(x))); \qquad  B_k= \Bbb
P_k(\mathfrak N_{R,\xi}(\theta_2(x)))  .
$$ Or
$$
 |\tau_k(A^p)-\tau_\omega(\theta_1(x))|<\frac 1 m; \quad |\tau_k(B^p)-\tau_\omega(\theta_2(x))|<\frac 1 m
$$ for all $1\le p\le m$ and $k\ge k_m.$
From Lemma 2, it follows that when $k$ is big enough  there is some
$U_k$ in $\mathcal U(k)$ such that
$$
\|U_kA_kU_k^*-B_k\|_2\le \delta.
$$
This implies that the $\delta$-orbit-covering number $\nu(\Bbb
P_k(\mathfrak N_{R,\xi}(x)),\delta)=1$. Therefore, $\mathcal
H_0^\omega(x) = 0$. By Theorem 1, we obtain that $\mathcal
H_0^\omega(\mathcal A)=0.$

\end{proof}

\begin{remark}
By the remark after Lemma 2, the following result also holds. If
$\mathcal N$ is a hyperfinite von Neumann algebra with a tracial
state, then $\mathcal H_0^\omega(\mathcal N)=0$.
\end{remark}

The following proposition is   Theorem 2.7 of \cite{V1}, whose proof
depends on the powerful tools from random matrices. An elementary
proof which is   based on the basic facts of unitary matrices can be
found in \cite{DH}.
\begin{proposition}  Let $L(F_n)$ be the free group factor  on $n$ generators with the tracial state $\tau$, and $u_1,\ldots, u_n$ be
the standard generators of $L(F_n)$.  For each $m, k\ge 1$, let
$$\begin{aligned}
\Omega_m(k) =& \{ (U_1, \ldots, U_n)\in \mathcal U(k)^n
 \ | \ |\tau_k(U_{i_1}^{\epsilon_1}\cdots U_{i_p}^{\epsilon_p})
  - \tau(u_{i_1}^{\epsilon_1}\cdots u_{i_p}^{\epsilon_p}) |<
 \frac 1 m \\
  &\qquad  \left . \text{ for all } \quad 1\le p\le m, 1\le i_1,\ldots, i_p\le n, \ \{\epsilon_1,\ldots, \epsilon_p\}\subset\{1,*\}\right \}.\end{aligned}
$$Then
$$
lim_{k\rightarrow \infty}\mu_k(\Omega_m(k))=1,
$$ where $\mu_k$ is normalized Haar measure on the compact group  $\mathcal U(k)^n$.

\end{proposition}

\begin{theorem}
Suppose $L(F_n)$ is the free group factor on $n$ generators with
$n\ge 2$. Then Voiculescu's embedding dimension
$\delta_0^\omega(L(F_n))\ge 2$ and $\mathcal
H_1^\omega(L(F_n))=\infty.$
\end{theorem}
\begin{proof}
 It follows from Proposition 1 that, for every $m\ge 1$, there are
 some positive integer $k_m$ and a sequence of subsets
 $\{\Omega_m(k)\}_{k=k_m}^\infty$ such that
$$  \mu_k(\Omega_m(k))\ge \frac 1 2, \qquad \quad  \text{for $k\ge k_m$},$$   where
$\mu_k$ is normalized Haar measure on the compact group $\mathcal
U(k)^n$.

Let $\xi=\{k_m\}_{m=1}^\infty$. It is easy to see that $\xi\in \Xi$.
For each $R>1$ and such $\xi$, consider the sequence
$\{\Sigma_k\}_{k=1}^\infty$ such that
  $$
     \Sigma_k= \Omega_m(k), \qquad \text { when } \ k_m\le k<
     k_{m+1}.
  $$
 It is not hard to verify that
 $$
   \prod_{k=1}^\infty \Sigma_k \subset \mathfrak N_{R,\xi}(
u_1,\ldots,u_n).
 $$ So,
$$
\Omega_m(k)=\Sigma_k \subset \Bbb P_k(\mathfrak N_{R,\xi}(
u_1,\ldots,u_n)) \qquad \text { when } \  k_m\le k< k_{m+1}.
$$Hence
$$
\mu_k(\Bbb P_k(\mathfrak N_{R,\xi}( u_1,\ldots,u_n)))\ge \frac 1 2
\qquad \text { for all } \ k\ge 1.
$$
Note there exists  constant   $c$  (not depending
  on   $k$) such that
$$\begin{aligned}
  \mu_k(Ball((U_1,\ldots,U_n),\delta)) \le  \left (  c\delta \right )^{ nk^2},
  \qquad   \forall \ 0<\delta<1 ,\end{aligned}
  $$ where  $Ball((U_1,\ldots,U_n),\delta)$ is a ball centered at $(U_1,\ldots,U_n)$
  with radius $\delta$ (with respect to $2$-norm) in $(\mathcal
  U(k))^n$.
 A standard argument on covering numbers shows that
$$
\delta_0^\omega(u_1,\ldots,u_n) =n.
$$ Similarly,
there exist a   constants    $ C$ (not depending
  on   $k$) such that
$$\begin{aligned}
      \mu_k(\mathcal U((U_1,\ldots,U_n),\delta)) \le  \left (  C\delta \right )^{ (n-1)k^2},
  \qquad   \forall \ 0<\delta<1 ,\end{aligned}
  $$ where   $\mathcal U((U_1,\ldots,U_n),\delta)$ is a unitary orbit  centered at $(U_1,\ldots,U_n)$
  with radius $\delta$ (with respect to $2$-norm)  in $(\mathcal
  U(k))^n$. A standard
  arguments on unitary orbit covering number shows
$$
\mathcal H_1^\omega(u_1,\ldots,u_n) =\infty.
$$ Thus, from Theorem 1, we have $\delta_0^\omega(L(F_n))\ge 2$ and $\mathcal
H_1^\omega(L(F_n))=\infty.$
\end{proof}

\begin{remark}
It seems that $\mathcal H_1^\omega$ does not provide us with more
insights into the isomorphism problem of free group factors because
$\mathcal H_1^\omega(L(F_n))=\infty$ for all $n\ge 2$. But it will
provide us with   useful information when von Neumann algebra is not
free group factors, which we will see in next sections.

\end{remark}

\section{Some lemmas on covering number and orbit covering number }
In this section, we are going to compute the covering numbers and
orbit-covering numbers of some sets.  We start with a  definition,
which is just for our convenience.
\begin{definition}
A unitary matrix $U$ in $\mathcal{M}_{k}(\mathbb{C})$ is a
\emph{Haar unitary matrix} if $\tau_{k}(U^{m})=0$ for all $1\leq
m<k$ and $\tau_{k}(U^{k})=1$.
\end{definition}

We have the following lemma.

\begin{lemma}
Let $V_{1} $ be a Haar unitary matrix  and $V_2$ be a unitary matrix in $\mathcal{M}_{k}%
(\mathbb{C})$. For every $\delta>0$, let
\[
\Omega(V_{1},V_{2};\delta)=\{U\in\mathcal{U}(k)\ |\ \Vert UV_{1}-V_{2}%
U\Vert_{2}\leq\delta\}.
\]
Then, for every $r>\delta $, there exists a set $\{Ball
(U_{\lambda}; \frac {4\delta} r)\}_{\lambda \in\Lambda}$ of
$\frac{4\delta}{r}$-balls in $\mathcal{U}(k)$ that cover
$\Omega(V_{1},V_{2};\delta)$ with the cardinality of $\Lambda$
satisfying $|\Lambda|\leq\left( \frac{3r}{2\delta}\right)
^{4rk^{2}}$.
\end{lemma}

\begin{proof}[Sketch of Proof]
\textbf{\ } Let $D$ be a diagonal unitary matrix, $diag(\lambda
_{1},\ldots
,\lambda _{k})$, where $\lambda _{j}$ is the $j$-th root of unity $1$. Since $%
V_{1} $ is a Haar unitary matrices, there exists $W_{1}$ in $%
\mathcal{U}(k)$ such that $V_{1}=W_{1}DW_{1}^{\ast }$. Assume that
$\mu_1,\ldots,\mu_k$ are the eigenvalues of $V_2$. Then there is
some unitary matrix $W_2$ such that $V_2= W_2 D_2 W_2^*$, where
$D_2=diag(\mu_{1},\ldots ,\mu_{k})$.
 Let $\tilde{\Omega}(\delta )=\{U\in \mathcal{U}%
(k)\ |\ \Vert UD-D_2U\Vert _{2}\leq \delta \}.$ Clearly  $\Omega
(V_{1},V_{2};\delta )=\{W_{2}^{\ast }UW_{1}|U\in \tilde{\Omega}(\delta)\}$; whence $%
\tilde{\Omega}(\delta)$ and $\Omega (V_{1},V_{2};\delta )$ have the
same covering numbers.

Let $\{e_{st}\}_{s,t=1}^{k}$ be the canonical system of matrix units
of $ \mathcal{M}_{k}(\mathbb{C})$. Let
\[
\begin{aligned}
\mathcal  S_1 = span \{e_{st} \  | \  |\lambda_s-\mu_t|< r \} \qquad
\mathcal S_2=M_k(\Bbb C) \ominus S_1.
\end{aligned}
\]%
For every $U=\sum_{s,t=1}^{k}x_{st}e_{st}$ in $\tilde{\Omega}(\delta
)$,
with $x_{st}\in \mathbb{C}$, let $T_{1}=\sum_{e_{st}\in \mathcal{S}%
_{1}}x_{st}e_{st}\in \mathcal{S}_{1}$ and $T_{1}=\sum_{e_{st}\in \mathcal{S}%
_{2}}x_{st}e_{st}\in \mathcal{S}_{2}$. But
\[\begin{aligned}
\delta ^{2}&\geq \Vert UD-D_2U\Vert
_{2}^{2}=\sum_{s,t=1}^{k}|(\lambda _{s}-\mu_{t})x_{st}|^{2}\geq
\sum_{e_{st}\in \mathcal{S}_{2}}|(\lambda
_{s}-\mu_{t})x_{st}|^{2}\\ &\geq r^{2}\sum_{e_{st}\in \mathcal{S}%
_{2}}|x_{st}|^{2}=r^{2}\Vert T_{2}\Vert _{2}^{2}.
\end{aligned}\]%
Hence $\Vert T_{2}\Vert _{2}\leq \frac{\delta }{r}$. Note that
$\Vert
T_{1}\Vert _{2}\leq \Vert U\Vert _{2}=1$ and $dim_{\mathbb{R}}{}\mathcal{S}%
_{1}\leq 4rk^{2} $ (see \cite{GS2}). By standard arguments on
covering numbers, we know
that $\tilde{\Omega}(\delta)$ can be covered by a set   $\{Ball(A^{\lambda };\frac{%
2\delta }{r})\}_{\lambda \in \Lambda }$ of $\frac{2\delta }{r}$-balls in $ \mathcal{M}_{k}(\mathbb{C})$ with $|\Lambda |\leq \left( \frac{3r}{2\delta }%
\right) ^{4rk^{2}}.$ Because $\tilde{\Omega}(\delta )\subset \mathcal{U}%
(k)$, after replacing   $A^{\lambda }$ by a unitary $U^{\lambda }$ in $%
Ball(A^{\lambda },\frac{2\delta }{r})$, we obtain a set $\{Ball
(U_{\lambda };\frac{4\delta }{r} )\}_{\lambda \in \Lambda }$ of
$\frac{4\delta }{r}$-balls in $\mathcal{U}(k)$ that cover
$\tilde{\Omega}(\delta )$ with the cardinality of $\Lambda $
satisfying $|\Lambda |\leq \left( \frac{3r}{2\delta }\right) ^{4r
k^{2}}$. Therefore the same result holds for $\Omega
(V_{1},V_{2};\delta )$.
\end{proof}

With the notations as above, we have following lemmas.
\begin{lemma} Suppose $R>1$, $0< r, \delta<1$. Let $r_1=\frac
{r\delta }{128R}.$ Suppose $\Gamma$ is a subset of $$(\mathcal
M_k(\Bbb C)^n\times \mathcal M_k(\Bbb C)\times \mathcal M_k(\Bbb
C)\times \mathcal M_k(\Bbb C))_R$$ such that   every
$$
(A_1,\ldots,A_n, U, V, W) \in \Gamma \subset (\mathcal M_k(\Bbb
C)^n\times \mathcal M_k(\Bbb C)\times \mathcal M_k(\Bbb C)\times
\mathcal M_k(\Bbb C))_R
$$ satisfies
\begin{enumerate}
\item [(i)] There exist some Haar unitary matrix $V_1$ and unitary matrices $U_1$, $W_1$ in $\mathcal U(k)$
such that $\|V-V_1\|_2<r_1$, $\|U-U_1\|_2< r_1$ and
$\|W-W_1\|_2<r_1;$
\item [(ii)] $\|UV-WU\|_2<r_1.$
\end{enumerate}
Let $$
\begin{aligned}
\Gamma_1 & = \{(A_1,\ldots,A_n, U) \in  (\mathcal M_k(\Bbb
C)^n\times \mathcal M_k(\Bbb C)) \ | \\ &\qquad \qquad\qquad\
\exists \ V, W \text { such
that } (A_1,\ldots,A_n, U, V, W) \in \Gamma \}\\
\Gamma_2 &= \{(A_1,\ldots,A_n,V, W) \in  (\mathcal M_k(\Bbb
C)^n\times \mathcal M_k(\Bbb C)\times \mathcal M_k(\Bbb C)) \ | \\
&\qquad \qquad \qquad  \ \exists \ U,
  \text { such that } (A_1,\ldots,A_n, U, V, W) \in \Gamma \}.
\end{aligned}
$$ Then  we have
$$
\nu(\Gamma_1,\delta) \le \nu(\Gamma_2,  \frac{r\delta}{128} ) \cdot
\left ( \frac {24 }{\delta}\right)^{4rk^2},
$$ where $\nu(\Gamma_1,\delta)$, or $\nu(\Gamma_2, \frac{r\delta}{128} )$, is the unitary orbit covering number of
the set $\Gamma_1$, or $\Gamma_2$ respectively, with radius
$\delta_1$, or $  \frac{r\delta}{128} $ respectively.

\end{lemma}
\begin{proof}By the definition of   unitary orbit covering number,
we know that $\Gamma_2$ can be covered by a collection of
$\frac{r\delta}{128}$-orbit-balls $\{\mathcal U(A_1^\lambda, \ldots,
A_n^\lambda, V^\lambda,W^\lambda; \frac{r\delta}{128})
\}_{\lambda\in\Lambda}$ such that the cardinality of $\Lambda$
satisfies $   | \Lambda|= \nu(\Gamma_2,\frac{r\delta}{128}). $

For every $ (A_1,\ldots,A_n, U, V, W) \in \Gamma$, we know
\begin{enumerate}
  \item [(i)] $(A_1,\ldots, A_n, V,W)$ is contained in $\Gamma_2$;
  \item [(ii)] There exist some Haar unitary matrix $V_1$ and unitary matrices $U_1$, $W_1$ in $\mathcal U(k)$
such that $\|V-V_1\|_2<r_1$, $\|U-U_1\|_2< r_1$ and
$\|W-W_1\|_2<r_1;$
\item [(iii)] $\|UV-WU\|_2<r_1.$
\end{enumerate}
From (i), it follows that there are some $\lambda$ in $\Lambda$ such
that
$$
(A_1,\ldots, A_n, V,W) \in \mathcal U(A_1^\lambda, \ldots,
A_n^\lambda, V^\lambda,W^\lambda;\frac{r\delta}{128}).
$$ Or, there is some
 unitary matrix $X$ in $\mathcal U(k)$ such that
$$
\| (A_1,\ldots, A_n, V,W)- X (A_1^\lambda, \ldots, A_n^\lambda,
V^\lambda,W^\lambda)X^*\|\le \frac{r\delta}{128}.
$$
Combining with (ii) and (iii), we obtain that
$$ \|U_1XV^\lambda X^* -XW^\lambda X^*U_1 \|_2< r_1 + 2 \frac{r\delta}{128} + 2R r_1,     $$
and
$$ \|V_1- XV^\lambda X^*\|_2\le r_1+  \frac{r\delta}{128}, \qquad \|W_1- XW^\lambda X^*\|_2\le r_1+  \frac{r\delta}{128} .   $$
It follows that there are some Haar unitary matrix $\tilde
V^\lambda$ and unitary matrix $\tilde W^\lambda$ in $\mathcal U(k)$
such that
$$\|\tilde V^\lambda- V^\lambda\|_2 \le r_1+ \frac{r\delta}{128}, \qquad \|\tilde W^\lambda- W^\lambda\|_2 \le r_1+ \frac{r\delta}{128}.$$
Replace this $V^\lambda$ by such Haar unitary matrix $\tilde
V^\lambda$ and this $W^\lambda$ by such unitary matrix $\tilde W$,
when $V^\lambda$ is not a Haar unitary matrix and $W^\lambda$ is not
a unitary matrix. Therefore, we have
$$ \|U_1X\tilde V^\lambda X^* -X\tilde W^\lambda X^*U_1 \|_2< 2(r_1+  \frac{r\delta}{128})+
(r_1 + 2 \frac{r\delta}{128} + 2R r_1)\le   \frac{r\delta}{16}.
$$
 On the other hand, it follows from Lemma 2, there exists a
set $\{Ball (U^{\lambda\sigma};\frac \delta 4)\}_{\sigma\in\Sigma}$
of $\frac \delta 4$-balls that covers $\Omega(\tilde V^\lambda,
\tilde W^\lambda; \frac{r\delta}{16})$ such that $ |\Sigma|\le\left
( \frac {24 }{\delta}\right)^{4rk^2}. $ Thus there is some $\sigma$
in $\Sigma$ such that
$$\|X^*U_1X-U^{\lambda\sigma}\|\le \frac \delta 4.$$ It
induces that
$$
 \| (A_1,\ldots,A_n, U, V, W) - X(A_1^\lambda, \ldots,
 A_n^\lambda,U^{\lambda\sigma},
V^\lambda,W^\lambda) X^*\|_2 \le \frac{r\delta}{128} + r_1 +\frac
\delta 4\le \frac \delta 2.
$$ From the definition of unitary orbit covering number it
follows that
$$
(A_1,\ldots,A_n, U, V, W)\in U(A_1^\lambda, \ldots,
 A_n^\lambda,U^{\lambda\sigma},
V^\lambda,W^\lambda; \delta)
$$ Hence,
$$
\nu(\Gamma_1, \delta)\le \nu(\Gamma,\delta) \le |\Lambda| |\Sigma|
\le \nu(\Gamma_2, \frac{r\delta}{128}) \cdot \left ( \frac {24
}{\delta}\right)^{4rk^2}.
$$
\end{proof}

\begin{lemma} Suppose $R>1$, $0< r, \delta<1$. Let $r_1=\frac
{r\delta }{96R}.$ Suppose $\Gamma$ is a subset of $$(\mathcal
M_k(\Bbb C)^n\times (\mathcal M_k(\Bbb C)^p\times \mathcal M_k(\Bbb
C))_R$$ such that   every
$$
(A_1,\ldots,A_n, C_1,\ldots,C_p, U) \in \Gamma \subset (\mathcal
M_k(\Bbb C)^n\times \mathcal M_k(\Bbb C)\times \mathcal M_k(\Bbb
C)\times \mathcal M_k(\Bbb C))_R
$$ satisfies that there exists some Haar unitary matrix $U_1$ in
$\mathcal U(k)$ such that   $$\|U-U_1\|_2< r_1.$$  Let $$
\begin{aligned}
\Gamma_1 & = \{(A_1,\ldots,A_n, U) \in  (\mathcal M_k(\Bbb
C)^n\times \mathcal M_k(\Bbb C)) \ | \\ &\qquad \qquad\qquad\
\exists \ C_1, \ldots, C_p \text { such
that } (A_1,\ldots,A_n, C_1,\ldots,C_p,U) \in \Gamma \}\\
\Gamma_2 &= \{(C_1,\ldots, C_p,U) \in  (\mathcal M_k(\Bbb C)^p\times
\mathcal M_k(\Bbb C)) \ | \\ &\qquad \qquad \qquad  \ \exists \
A_1,\ldots, A_n
  \text { such that } (A_1,\ldots,A_n, C_1,\ldots, C_p,U) \in \Gamma \}.
\end{aligned}
$$ Then  we have
$$
\nu(\Gamma ,2n\delta) \le \nu(\Gamma_1,  \frac{r\delta}{96} )\cdot
\nu(\Gamma_2,  \frac{r\delta}{96} ) \cdot \left(
\frac{18R}{\delta}\right)  ^{4rk^{2}}
$$ where $\nu(\Gamma ,2n\delta)$, or $\nu(\Gamma_i, \frac{r\delta}{96} )$  $(i=1,2)$, is the unitary orbit covering number of
the set $\Gamma $, or $\Gamma_i$  $(i=1,2)$  respectively, with
radius $2n\delta_1$, or $  \frac{r\delta}{96} $ respectively.

\end{lemma}

\begin{proof}
By the definition of unitary orbit covering numbers of the sets
$\Gamma_1$ and $\Gamma_2$,
 there exists a set $\{\mathcal{U}%
(B_{1}^{\lambda},\ldots,B_{n}^{\lambda},U_{\lambda};\frac{r\delta}%
{96R})\}_{\lambda\in\Lambda }$ of $\frac{r\delta}{96R}$-orbit-balls
in $(\mathcal{M}_{k}(\mathbb{C})^{n+1})_R$   covering $\Gamma_1$
with $|\Lambda |=\nu(\Gamma_1,\frac{r\delta}{96R})$.  Also there
exists a set $\{\mathcal{U}(D_{1}^{\sigma},\ldots,D_{p}^{\sigma
},U_{\sigma};\frac{r\delta}{96R})\}_{\sigma\in\Sigma }$ of
$\frac{r\delta }{96R}$-orbit-balls in
$(\mathcal{M}_{k}(\mathbb{C})^{p+1})_R$ that cover $\Gamma_2$ with
$|\Sigma |=\nu (\Gamma_2,\frac{r\delta}{96R})$.

For each $(A_{1},\ldots,A_{n},C_{1},\ldots,C_{p},U)$ in $\Gamma$, we
know the following hold. \begin{enumerate} \item [(i)]
$(A_{1},\ldots,A_{n},U)$ is contained in $\Gamma_1$;
\item [(ii)] $(C_{1},\ldots,C_{p},U)$ is contained in $\Gamma_2$;
\item [(iii)] There exists some Haar unitary matrix $U_1$ in $\mathcal U(k)$
such that \\ \qquad \qquad \makebox[2cm]{}  \qquad $\| U-U_1\|_2\le
r_1.$
\end{enumerate}
From (i) and (ii), there exist some $\lambda$ in $\Lambda$, $\sigma$
in $\Sigma$ such that
$$
\begin{aligned}
  (A_{1},\ldots,A_{n},U) &\in U(B_1^\lambda, \ldots,
  B_n^\lambda,U_\lambda; \frac {r\delta}{96R})\\
  (C_{1},\ldots,C_{p},U) &\in U(D_1^\sigma, \ldots,
  D_n^\sigma,U_\sigma; \frac {r\delta}{96R})
\end{aligned}
$$
i.e., there exist unitary matrices $W_1, W_2$ in $\mathcal U(k)$
such that
\[
\begin{aligned}
&\|(A_1,\ldots,A_n,U) - (W_1B_1^\lambda W_1^*, \ldots,
W_1B_n^\lambda
W_1^*, W_1U_\lambda W_1^*)\|_2 \le \frac {r\delta}{96R}\\
& \|( C_1,\ldots,C_p, U)-(W_2D_1^\sigma W_2^*, \ldots, W_2D_p^\sigma
W_2^*, W_2U_\sigma W_2^*)\|_2\le \frac {r\delta}{96R}.
\end{aligned}
\]
Combining with (iii), we have that
$$
\|U_1-W_1U_\lambda W_1^*\|_2\le r_1 + \frac {r\delta}{96R}\le \frac
{r\delta}{48R}, \quad \|U_1-W_2U_\sigma W_2^*\|_2\le r_1 + \frac
{r\delta}{96R}\le \frac {r\delta}{48R}.
$$
Therefore there are  two Haar unitary matrices $V_\lambda, V_\sigma$
such that
$$
\| U_\lambda - V_\lambda\|_2\ \le \frac {r\delta}{48R}, \qquad \|
U_\sigma- V_\sigma\|_2 \le \frac {r\delta}{48R}.
$$
Replace these $U_\lambda, U_\sigma$ by Haar unitary matrices
$V_\lambda, V_\sigma$ when $U_\lambda, U_\sigma$ are not Haar
unitary matrices.
 It follows
\[
\Vert
W_{2}^{\ast}W_{1}V_{\lambda}-V_{\sigma}W_{2}^{\ast}W_{1}\Vert_{2}=\Vert
W_{1}V_{\lambda}W_{1}^{\ast}-W_{2}V_{\sigma}W_{2}^{\ast}\Vert_{2}\leq
\frac{r\delta}{12R}.
\]
Since $V_{\lambda},V_{\sigma}$ are Haar unitary matrices in
$\mathcal{M}_{k}(\mathbb{C})$, by Lemma 3 we know that there exists
a set $\{Ball (U_{\lambda\sigma\gamma}; \frac{\delta}
{3R})\}_{\gamma\in\mathcal{I}_{k}}$ of $\frac{\delta}  {3R}$-balls
in $\mathcal{U}(k)$ that cover $\Omega(V_{\lambda},V_{\sigma
};\frac{r\delta}{12R})$ with the cardinality of $\mathcal{I}_{k}$
never exceeding $\left( \frac{18R}{\delta}\right)  ^{4rk^{2}}.$ Then
there exists some $\gamma\in\mathcal{I}_{k}$ such that $\Vert
W_{2}^{\ast}W_{1}-U_{\lambda\sigma\gamma}\Vert_{2}\leq\frac{\delta}{3R}$.
This in turn implies
\[
\begin{aligned}
\|(A_1,\ldots, A_n,C_1,\ldots,C_p, U) - &
(W_2U_{\lambda\sigma\gamma}B_1^\lambda
U_{\lambda\sigma\gamma}^*W_2^*, \ldots,
W_2U_{\lambda\sigma\gamma}B_n^\lambda
U_{\lambda\sigma\gamma}^*W_2^*, \\
&  \quad\quad   W_2D_1^\sigma W_2^*, \ldots, W_2D_p^\sigma
W_2^*,W_2U_\sigma W_2^* )\|_2\le  n\delta\end{aligned}
\]
for some
$\lambda\in\Lambda_{k},\sigma\in\Sigma_{k},\gamma\in\mathcal{I}_{k}$
and $W_{2}\in\mathcal{U}(k)$, i.e.,
\[
(A_{1},\ldots,A_{n},C_{1},\ldots,C_{p},U)\in\mathcal{U}(U_{\lambda\sigma
\gamma}B_{1}^{\lambda}U_{\lambda\sigma\gamma}^{\ast},\ldots,U_{\lambda
\sigma\gamma}B_{n}^{\lambda}U_{\lambda\sigma\gamma}^{\ast},D_{1}^{\sigma
},\ldots,D_{p}^{\sigma},U_{\sigma};2n\delta).
\]
From the definition of unitary orbit covering number it follows that
$$
\nu(\Gamma, \delta) \le |\Lambda|\cdot |\Sigma| \cdot \left(
\frac{18R}{\delta}\right)  ^{4rk^{2}} =\nu(\Gamma_1,
\frac{r\delta}{96} )\cdot \nu(\Gamma_2,  \frac{r\delta}{96} ) \cdot
\left( \frac{18R}{\delta}\right)  ^{4rk^{2}}
$$
\end{proof}

\section{The computation of $\mathcal H_1^\omega$ for some finite von Neumann algebras}

In this section, we are going to compute $\mathcal
H_1^\omega(\mathcal N)$ for a finite von Neumann algebra $\mathcal
N$, including type II$_1$ factors with property $\Gamma$, with
Cartan subalgebras and nonprime type II$_1$ factors.

\subsection{Embedding dimension}

The same strategy as in the proof of Theorem 1 can be used to prove
the following theorem whose proof is skipped here.
\begin{theorem}
Suppose $\mathcal N$ is a finitely generated von Neumann algebra
with a tracial state $\tau$ and $\mathcal N$ can be faithfully
trace-preserving embedded into $\mathcal M_k(\Bbb C)^\omega$.
Suppose $\{\mathcal N_j\}_{j=1}^\infty$ is a increasing sequence of
von Neumann subalgebras of $\mathcal N$ such that $\mathcal
N=\overline{\cup_{j=1}^\infty \mathcal N_j}^{SOT}$. Then, for each
$s\ge 0$,
$$0\le \ \mathcal H_s^\omega(\mathcal N) \ \le \ \liminf_{j\rightarrow \infty}
\mathcal H_s^\omega(\mathcal N_j).$$
\end{theorem}

\begin{definition}
Suppose that $\mathcal{N}$ is a diffuse von Neumann algebra with a
tracial state $\tau$. Then a unitary element $u$ in $\mathcal{N}$ is
called a \emph{Haar unitary} if $\tau(u^{m})=0$ when $m\neq0$.
\end{definition}

\begin{theorem}
Suppose $\mathcal{N}$ is a diffuse finitely generated von Neumann
algebra with a tracial state $\tau$ and $\mathcal N$ can be
faithfully trace-preserving embedded into $\mathcal M_k(\Bbb
C)^\omega$. Suppose $\mathcal{N}_1$ is a diffuse von Neumann
subalgebra of $\mathcal{N}$ and $u$ is a unitary element in
$\mathcal{N}$ such that  $\{\mathcal{N}_1,u\}$ generates
$\mathcal{M}$ as a von Neumann algebra. If there exist Haar unitary
elements $v_{1},v_{2},\ldots$ and unitary elements
$w_{1},w_{2},\ldots$ in
$\mathcal{N}_1$ such that $\left\Vert v_{n}u-uw_{n} \right\Vert _{2}%
\rightarrow 0$, then $\mathcal H_{1}^\omega( \mathcal{N} )\le
\mathcal H_1^\omega(\mathcal N_1)$. In particular, if there are Haar
unitary elements $v,w$ in $\mathcal{N}_1,$ such that $vu=uw$, then
$\mathcal H_1^\omega\left( \mathcal{N}\right) \le \mathcal
H_1^\omega(\mathcal N_1).$
\end{theorem}

\begin{proof}
\textbf{\ }
Suppose that $\{x_{1},\ldots,x_{n}\}$ is a family of generators of
$\mathcal{N}_{1}$. Then we know that
$\{x_{1},\ldots,x_{n}, u\}$ is a family of generators of
$\mathcal{M}$.

For every $0<\delta<1$, $0<r<1$, there exist an integer $p>0$ and a
Haar unitary element  $v_{p}$, and a unitary element $w_{p}$ in
$\mathcal{N}_1$ such that
\[
\|v_{p}u-uw_{p}\|_{2}< \frac{r\delta}{130}.
\]
Note that $\{x_{1},\ldots,x_{n}, v_{p},w_{p}\} $ is also a family of
generators of $\mathcal{N}_1$.

For $R>1 $, $\xi=\{k_m\}_{m=1}^\infty$ in $\Xi$, let  $r_1=\frac
{r\delta }{128R}.$  By the definition of $$\mathfrak
N_{R,\xi}(x_1,\ldots,  x_n,u,v_p,w_p)$$ and Lemma 2, there is some
$m_0\ge 0$ such that every
$$
((A_{1,k})_k,\ldots, (A_{n,k})_k, (U_{k})_k, (V_{k})_k,
(W_{k})_k)\in \mathfrak N_{R,\xi}(x_1,\ldots, x_n,u,v_p,w_p),
$$ satisfies, when $k\ge k_{m_0}$,
\begin{enumerate}
\item [(i)] there exist some Haar unitary matrix $V_{1,k}$ and unitary matrices $U_{1,k}$, $W_{1,k}$ in $\mathcal U(k)$
such that$$\text{ $\|V_k-V_{1,k}\|_2<r_1$, $\|U_k-U_{1,k}\|_2< r_1$
and $\|W-W_{1,k}\|_2<r_1;$} $$
\item [(ii)] $\|U_kV_k-W_kU_k\|_2<r_1.$
\end{enumerate}
Apply Lemma 3 by letting $\Gamma =\Bbb P_k(\mathfrak
N_{R,\xi}(x_1,\ldots, x_n,u,v_p,w_p))$. We get that
$$\begin{aligned}
\nu(\Bbb P_k(\mathfrak N_{R,\xi}&(x_1,\ldots,x_n,u:
v_p,w_p)),\delta) \\ & \le \nu(\Bbb P_k(\mathfrak N_{R,\xi}( x_1
,\ldots,  x_n, v_p,w_p:u)), \frac{r\delta}{128} ) \cdot \left (
\frac {24 }{\delta}\right)^{4rk^2}
\\ & \le \nu(\Bbb P_k(\mathfrak
N_{R,\xi}(x_1,\ldots,x_n, v_p,w_p)), \frac{r\delta}{128} ) \cdot
\left ( \frac {24 }{\delta}\right)^{4rk^2}.\end{aligned}
$$
 Hence, by the definition of   embedding dimension, we have shown
\[
\begin{aligned}
0\le  &\mathcal H_1^\omega (x_1, \ldots,x_n, u:
v_p,w_p;R,\xi,\delta)\\
&\le  \lim_{k\rightarrow \omega} \frac {\log(\nu(\Bbb P_k(\mathfrak
N_{R,\xi}(x_1,\ldots,x_n,  v_p,w_p)), \frac{r\delta}{128} ) \cdot
\left ( \frac {24 }{\delta}\right)^{4rk^2})}{ k^2 }\\
&= \lim_{k\rightarrow \omega} \frac {\log(\nu(\Bbb P_k(\mathfrak
N_{R,\xi}(x_1,\ldots,x_n, v_p,w_p)), \frac{r\delta}{128} ) )}{k^2 }+
4r \cdot  ({\log 24-\log\delta})
\\&\le \mathcal H_s^\omega(x_1,\ldots,x_n, v_p,w_p)+ 4r \cdot  ({\log
24-\log\delta})\\
&= \mathcal H_s^\omega(x_1,\ldots,x_n)+ 4r \cdot  ({\log
24-\log\delta}).
\end{aligned}
\]
Because $r$ is an arbitrarily small positive number, we have
$$\mathcal H_1^\omega (x_1, \ldots,x_n, u: v_p,w_p)\le \mathcal
H_s^\omega(x_1,\ldots,x_n).$$ By Lemma 1 and Theorem 1, $\mathcal
H_1^\omega(\mathcal N)\le  \mathcal H_1^\omega(\mathcal N_1)$.
\end{proof}

\begin{theorem}
Suppose $\mathcal{N}$ is a finitely generated von Neumann algebra
with a tracial state $\tau$ and $\mathcal N$ can be faithfully
trace-preserving embedded into $\mathcal M_k(\Bbb C)^\omega$.
Suppose $\mathcal{N}$ is generated by von Neumann subalgebras $\mathcal{N}%
_{1}$ and $\mathcal{N}_{2}$ of $\mathcal{N}$. If   $\mathcal{N}%
_{1}\cap\mathcal{N}_{2}$ is a diffuse von Neumann subalgebra of
$\mathcal{N}$, then $$\mathcal H_1^\omega( \mathcal{N})\le \mathcal
H_1^\omega( \mathcal{N}_1)+ \mathcal H_1^\omega( \mathcal{N}_2) .$$
\end{theorem}

\begin{proof}
\textbf{\ } Suppose that $\{x_{1},\ldots,x_{n}\}$ is a family of
generators of $\mathcal{N}_{1}$ and $\{y_{1},\ldots,y_{p}\}$ a
family of generators of $\mathcal{N}_{2}$. Since
$\mathcal{N}_{1}\cap\mathcal{N}_{2}$ is a diffuse von
Neumann subalgebra, we can find a Haar unitary $u$ in $\mathcal{N}_{1}%
\cap\mathcal{N}_{2}$.

For every $R>1+\max_{1\leq i\leq n,1\leq j\leq p}\{\Vert
x_{i}\Vert,\Vert y_{j}\Vert\}$, $0<\delta<\frac{1}{2n}$, $0<r<1$ and
$\xi\in \Xi$, let    $r_1=\frac {r\delta }{96R}.$ From the
definition of  $ \mathcal N_{R,\xi}(x_1,\ldots,x_n,y_1,\ldots,y_p,u)
$ and Lemma 2, there is some $m_0\ge 0$ such that, every
$$
((A_{1,k})_k,\ldots,(A_{n,k})_k,
(B_{1,k})_k,\ldots,(B_{p,k})_k,(U_k)_k) \in  \mathcal
N_{R,\xi}(x_1,\ldots,x_n,y_1,\ldots,y_p,u)
$$ satisfies that, when $k\ge k_{m_0}$, there exists a Haar
unitary matrix $U_{1,k}$ such that $$\|U_k-U_{1,k}\|_2\le r_1,$$
 Apply Lemma 5 by letting $\Gamma= \Bbb P_k(\mathcal
N_{R,\xi}(x_1,\ldots,x_n,y_1,\ldots,y_p,u))$. We get that
$$\begin{aligned}
\nu( \Bbb P_k(\mathcal
N_{R,\xi}&(x_1,\ldots,x_n,y_1,\ldots,y_p,u)),2n\delta) \le \nu( \Bbb
P_k(\mathcal N_{R,\xi} (x_1,\ldots,x_n,u:
y_1,\ldots,y_p)),\frac{r\delta}{96})\\
& \quad \cdot \nu( \Bbb P_k(\mathcal N_{R,\xi} ( y_1,\ldots,y_p,u:
x_1,\ldots,x_n)),\frac{r\delta}{96}) \cdot \left( \frac{18R}{\delta}\right)  ^{4rk^{2}}\\
& \le \nu( \Bbb P_k(\mathcal N_{R,\xi} (x_1,\ldots,x_n,u
)),\frac{r\delta}{96}) \cdot \nu( \Bbb P_k(\mathcal N_{R,\xi} (
y_1,\ldots,y_p,u )),\frac{r\delta}{96})  \cdot \left(
\frac{18R}{\delta}\right)  ^{4rk^{2}}
\end{aligned}
$$
Now it not hard to show that
$$\begin{aligned}
\mathcal H_1^\omega&(x_1,\ldots,x_n,y_1,\ldots,y_p,u;
R,\xi,2n\delta)\\
& \le \mathcal H_1^\omega (x_1,\ldots, x_n,u)+ \mathcal H_1^\omega
(y_1,\ldots, y_p,u)+ 4r \cdot (\log 18R -\log \delta).\end{aligned}
$$
Because $r$ is an arbitrarily small positive number, we have
$$\begin{aligned}
\mathcal H_1^\omega (x_1,\ldots,x_n,y_1,\ldots,y_p,u) \le \mathcal
H_1^\omega&(x_1,\ldots, x_n,u)+ \mathcal H_1^\omega (y_1,\ldots,
y_p,u).\end{aligned}
$$ By Theorem 1, we obtain,
$$
\mathcal H_1^\omega(\mathcal N) \le \mathcal H_1^\omega(\mathcal
N_1) + \mathcal H_1^\omega(\mathcal N_2).
$$

\end{proof}

Now we are able to going to compute   values of $\mathcal
H_1^\omega$ for many specific type II$_1$ factors.

\begin{theorem}
Suppose that $\mathcal N$ is a type II$_1$ factor with Cartan
subalgebras and $\mathcal N$ can be faithfully trace-preserving
embedded into $\mathcal M_k(\Bbb C)^\omega$. Then $\mathcal
H_1^\omega(\mathcal N)=0$.
\end{theorem}
\begin{proof}It follows from \cite {Po3} that $\mathcal N$ is generated by two
self-adjoint elements of $\mathcal N$(see also \cite{Sh}). Note
there is a maximal abelian von Neumann subalgebra $\mathcal A$ of
$\mathcal N$ such that $N(\mathcal A)$, the normalizers of $\mathcal
A$ in $\mathcal N$, generates $\mathcal N$  where $N(\mathcal A)$ is
the group of all these unitary elements $u$ in $\mathcal N$ such
that $u^*\mathcal A u=\mathcal A.$ Now it follows directly from
Theorem 1, Theorem 4 and Theorem 5 that $\mathcal
H_1^\omega(\mathcal N)=0$.
\end{proof}

\begin{corollary}
Suppose that $\mathcal R$ is the hyperfinite type II$_1$ factor.
Then $\mathcal H_1^\omega(\mathcal R)=0$.
\end{corollary}

\begin{theorem}
Suppose that $\mathcal N$ is type II$_1$ factor with property
$\Gamma$ and $\mathcal N$ can be faithfully trace-preserving
embedded into $\mathcal M_k(\Bbb C)^\omega$. Then $\mathcal
H_1^\omega(\mathcal N)=0$.
\end{theorem}
\begin{proof}
It follows from \cite {GePo} that $\mathcal N$ is generated by two
self-adjoint elements (see also \cite{Sh}). From Theorem 5.3 of
\cite {ChSm}, it follows that there is a hypefinite II$_1$ subfactor
$\mathcal R$ of $\mathcal N$ such that $\mathcal N'\cap \mathcal
R^\omega$ is diffuse. Now it follows from Corollary 1, Theorem 4 and
Theorem 5 that $\mathcal H_1^\omega(\mathcal N)=0.$
\end{proof}

The following two theorems also follows directly from Theorem 1,
Theorem 4 and Theorem 5, whose proofs are skipped here.
\begin{theorem}
Suppose that $\mathcal N$ is a non-prime type II$_1$ factor, i.e.
the tensor product of two type II$_1$ subfactors, and $\mathcal N$
can be faithfully trace-preserving embedded into $\mathcal M_k(\Bbb
C)^\omega$. Then $\mathcal H_1^\omega(\mathcal N)=0$.
\end{theorem}
\begin{theorem}Suppose that $SL(2n+1,\Bbb Z)$ is the special linear
group with the integer entries. Then
 $\mathcal H_1^\omega(L(SL(2n+1,\Bbb Z)))=0$.
\end{theorem}

\begin{remark}
We know that $\mathcal H_s^\omega(\mathcal N)$ is a decreasing
function of $s\ge 0$. The following question is of interest to us.
Suppose that $\mathcal N$ is a type II$_1$ factor with $\mathcal
H_1^\omega(\mathcal N)=0$. Can we find some number $t$ such that
$\mathcal H_t^\omega(\mathcal N)>0$? How about $L(F_2)\otimes
L(F_2)$? How about type II$_1$ factors with Cartan subalgebras?
\end{remark}

\subsection{Voiculescu's embedding dimension  of thin factors}

 The concept of thin factor was introduced by S. Popa (also see \cite {GePo}). It was known in \cite{GePo} that free group factors on $n$
 generators with $n\ge 4$ are not thin factors. This concept was
 further generalized to $\mathfrak K$-thin in \cite{HaSh}. It was
 shown there that free group factors on $n$
 generators with $n\ge 4$ are not  $\mathfrak K$-thin factors. Here, we are going to consider Voiculescu's embedding dimensions of
 these ``thin" factors.

\begin{theorem}Suppose that $\mathcal N$ is a finitely generated type II$_1$
 factor with a tracial state $\tau$ and $\mathcal N$ can be
 faithfully trace-preserving embedded into $\mathcal M_k(\Bbb
 C)^\omega$. Suppose there exist two subalgebras $\mathcal N_0$,
 $\mathcal N_1$ and $n$-vectors $\xi _{1},\ldots,\xi_{n}$$)$ in $L^{2}\left(
\mathcal{M},\tau\right)  $ such that $0\le \mathcal
H_1^\omega(\mathcal N_0),\mathcal H_1^\omega(\mathcal N_1) <\infty$,
and $\overline{span}^{\Vert\cdot\Vert_{2}}\mathcal{N}_{0}\{\xi
_{1},\ldots,\xi_{n}\}\mathcal{N}_{1}=L^{2}\left(
\mathcal{M},\tau\right) )$.
  Then   Voiculescu's embedding dimension of $\mathcal N$ satisfies $$
\delta_0^\omega(\mathcal N)\le 1+ 2n
  $$
\end{theorem}

\begin{proof}The proof of this theorem is just a slight modification
of the proof of    Theorem 7 in \cite{HaSh}.

\end{proof}

The following corollary follows easily from the preceding theorem.
\begin{corollary}
Suppose $\mathcal N$ is a type II$_1$ factor with a simple masa or
is a thin factor and $\mathcal N$ can be faithfully embedded into
$\mathcal M_k(\Bbb C)^\omega$. Then Voiculescu's embedding dimension
$\delta_0(\mathcal N) \le 3.$
\end{corollary}

\end{document}